\documentclass[12pt]{amsart}
\usepackage{euscript,eufrak}  
\usepackage{graphicx,epic,eepic}  
\usepackage{amsfonts}
\usepackage{amssymb}
\usepackage{latexsym}
\usepackage{amscd}

\newtheorem{thm}{Theorem}[section]
\newtheorem{cor}[thm]{Corollary}

\newtheorem{lemma}[thm]{Lemma}
\newtheorem{prop}[thm]{Proposition}

\newtheorem{conj}[thm]{Conjecture}

\newenvironment{pf*}[1]{\proof[#1]}{\endproof}
\usepackage{euscript}

\usepackage[OT2,OT1]{fontenc}

\newcommand{\cal}[1]{{\mathcal #1}}

\renewcommand{\Bbb}[1]{{\mathbb #1}}

\theoremstyle{definition}

\newcommand\Prefix[3]{\vphantom{#3}#1#2#3}
\theoremstyle{remark}

\renewcommand{\mod}{\operatorname{mod}}
\newcommand{\tl}{\tilde}

\newcommand{\eps}{\epsilon}

\newcommand{\aaa}[1]{{{\mathbf{#1}}}}

\newcommand{\comm}[1]{}

\renewcommand{\Re}{\operatorname{Re}}
\renewcommand{\Im}{\operatorname{Im}}

\numberwithin{equation}{section}
\newcommand{\thmref}[1]{Theorem~\ref{#1}}

\newcommand{\cA}{{\cal A}}

\newcommand{\cU}{{\cal U}}

\newcommand{\cH}{{\cal H}}
\newcommand{\cR}{{\cal R}}
\newcommand{\cL}{{\cal L}}

\newcommand{\CC}{{\Bbb C}}

\newcommand{\ZZ}{{\Bbb Z}}
\newcommand{\NN}{{\Bbb N}}
\newcommand{\DD}{{\Bbb D}}

\newcommand{\cren}{\cR_{\text cyl}}

\newcommand{\field}[1]{\mathbb{#1}}

\begin{document}
\addtolength{\evensidemargin}{-0.7in}
\addtolength{\oddsidemargin}{-0.7in}

\title{Cylinder renormalization of Siegel disks.}
\author{Denis Gaidashev, Michael Yampolsky}
\thanks{The second author is partially supported by an NSERC operating grant}

\begin{abstract}
We study one of the central open questions in one-dimensional renormalization theory -- the
conjectural universality of golden-mean Siegel disks. We present an approach to the 
problem based on cylinder renormalization proposed by the second author. 
Numerical implementation of this approach relies on the Constructive Measurable Riemann 
Mapping Theorem proved by the first author. Our numerical study yields a convincing evidence
to support the Hyperbolicity Conjecture in this setting.
\end{abstract}

\date{\today}

\maketitle

\section{Introduction}

One of the central examples of universality on one-dimensional dynamics is
provided by Siegel disks of quadratic polynomials. Let us consider, for instance,
the mapping 
$$P_\theta(z)=z^2+e^{2\pi i\theta}z,\text{ where }\theta=(\sqrt{5}+1)/2$$
is the golden mean. By a classical result of Siegel, the dynamics of $P_\theta$
is linearizable near the origin. The Siegel disk of $P_\theta$, which we will further
denote $\Delta_\theta$ is the maximal neighborhood of zero in which a conformal
change of coordinates reduces $P_\theta$ to the form $w\mapsto e^{2\pi i\theta}w$.
By the results of Douady, Ghys, Herman, and Shishikura, the topological disk $\Delta_\theta$
extends up to the only critical point of $P_\theta$ and is bounded by  a Jordan curve. 

It has  been observed numerically (cf. the work of Manton and Nauenberg \cite{MN}), 
that the boundary of $\Delta_\theta$ is asymptotically
self-similar near the critical point. Moreover, the scaling factor is universal
in a large class of analytic mappings with a golden-mean Siegel disk.
In 1983 Widom \cite{Wi} defined a renormalization procedure for $P_\theta$ which
``blows up'' a part of the invariant curve $\partial \Delta_\theta$ near the critical
point, and conjectured that the renormalizations of $P_\theta$ converge to a fixed point.
In addition, he conjectured that in a suitable functional space this fixed point
is hyperbolic with one-dimensional unstable direction.

In 1986
MacKay and Persival \cite{MP} extended the conjecture to other rotation numbers,
postulating the 
existence of a hyperbolic renormalization horseshoe corresponding to Siegel disks of 
analytic maps,
analogous to the Lanford's horseshoe for critical circle maps \cite{La1,La2}. 

In 1994 Stirnemann \cite{Stir} gave a computer-assisted proof of the existence of a
renormalization fixed point with a golden-mean Siegel disk.
In 1998, McMullen \cite{McM} proved the asymptotic self-similarity of  golden-mean
Siegel disks in the quadratic family. He constructed a version of renormalization 
based on holomorphic commuting pairs of de~Faria \cite{dF1,dF2}, and
showed that the renormalizations of a quadratic polynomial with a golden Siegel disk
 near the critical point converge
to a fixed point geometrically fast. 
More generally, he constructed a renormalization horseshoe for bounded
type rotation numbers, and used renormalization to show that the Hausdorff dimension of
the corresponding quadratic Julia sets is strictly less than two.

Having thus attracted much attention, the hyperbolicity part of the conjecture of Widom for golden-mean
Siegel disks is still open.

In \cite{Ya1} the second author has introduced a new renormalization transformation $\cren$, which he called
the cylinder renormalization, and used it to prove the Lanford's Hyperbolicity
Conjecture for critical circle maps. The main advantage of $\cren$ over
the renormalization scheme based on commuting pairs is that this operator is analytic
in a Banach manifold of analytic maps of a subdomain of $\CC/\ZZ$. 
It is thus a natural setting to study the hyperbolic properties of a fixed point.
In the present paper we study the fixed point of the cylinder renormalization numerically, and empirically
confirm the hyperbolicity conjecture, as well as study the dynamical properties of the fixed point.
The main numerical challenge in working with cylinder renormalization is a change of coordinate
involved in its definition. It is defined implicitly, and uniformizes a dynamically defined
fundamental domain to the straight cylinder $\CC/\ZZ$. To handle it, we use the Constructive
Measurable Riemann Mapping Theorem developed for numerically solving the Beltrami partial differential equation
by the first author in \cite{Gay,GayKhmel}.

\section{Definition and main properties of the cylinder renormalization of Siegel disks.}
\label{section: cyl ren}

\noindent
{\bf Some functional spaces.} 
For a topological disk $W\subset\CC$ containing $0$ and $1$ we will denote 
${\aaa A}_W$ the Banach space of bounded analytic functions in $W$ equipped with the
sup norm. Let us denote ${\aaa C}_W$ the Banach subspace of 
${\aaa A}_W$ consisting of analytic mappings $h:W\to \CC$ such that
$h(0)=0$ and $h'(1)=0$. 

In the case when the domain $W$ is the disk $\DD_\rho$ of radius $\rho>1$ centered
at the origin, we will denote ${\bf A}_{\DD_\rho}\equiv {\bf A}_\rho$ and 
${\bf C}_{\DD_\rho}\equiv {\bf C}_\rho$.

For each $\rho>1$ we will also consider the 
collection ${\bf B}^1_\rho$ of analytic functions $f(z)$ defined on some neighborhood
of the origin with $f(0)=0$, equipped with the weighted $l_1$
norm on the coefficients of the Maclaurin's series:

 \begin{equation}
 \|f\|_\rho=\sum^{\infty}_{n=0}{\left| f^{(n)}(0) \right| \over n !} \rho^n.
 \end{equation}

\noindent
We will further denote ${\bf L}^1_\rho$ the subset of ${\bf B}^1_\rho$ consisting of maps $f$ with
the normalizing condition $f'(1)=0$.


The proof of the following elementary statement is left to the reader:

\begin{lemma}\label{Equivalence}
$ \ $  
\begin{itemize}
\item[1)]  Let $f \in {\bf L}^1_\rho$, then 
$\sup_{\field{D}_\rho} |f(z)| \le  \|f\|_\rho$;
\item[2)]  Let $f \in {\bf A}_{\rho'}$ and $\rho'>\rho$, then $ \|f\|_{\rho} \le {\rho \over \rho'-\rho}  \sup_{\field{D}_{\rho'}} |f(z)|.$
\end{itemize}
\end{lemma}

As an immediate consequence, we have:

\begin{cor}
${\bf L}_\rho^1$ is a Banach space.
\end{cor}

\noindent

\begin{figure}
  \begin{center}
    \setlength{\unitlength}{0.00083333in}
    \begingroup\makeatletter\ifx\SetFigFont\undefined%
    \gdef\SetFigFont#1#2#3#4#5{%
      \reset@font\fontsize{#1}{#2pt}%
      \fontfamily{#3}\fontseries{#4}\fontshape{#5}%
      \selectfont}%
    \fi\endgroup%
	{\renewcommand{\dashlinestretch}{30}
	  \begin{picture}(4930,1695)(0,-10)
	    \drawline(2022.000,445.000)(2064.009,379.174)(2116.489,321.349)
	    (2177.945,273.171)(2246.626,236.013)(2320.577,210.933)
	    (2397.693,198.644)(2475.777,199.497)(2552.606,213.467)
	    (2625.992,240.157)(2693.845,278.806)(2754.234,328.315)
	    (2805.438,387.272)(2846.000,454.000)
	    \drawline(2024.000,440.000)(2066.330,374.753)(2118.993,317.519)
	    (2180.499,269.916)(2249.110,233.289)(2322.888,208.673)
	    (2399.746,196.765)(2477.514,197.900)(2553.992,212.047)
	    (2627.020,238.806)(2694.533,277.420)(2754.623,326.798)
	    (2805.592,385.545)(2846.000,452.000)
	    \drawline(2018.000,1384.000)(2060.330,1318.753)(2112.993,1261.519)
	    (2174.499,1213.916)(2243.110,1177.289)(2316.888,1152.673)
	    (2393.746,1140.765)(2471.514,1141.900)(2547.992,1156.047)
	    (2621.020,1182.806)(2688.533,1221.420)(2748.623,1270.798)
	    (2799.592,1329.545)(2840.000,1396.000)
	    \drawline(2841.000,1396.000)(2807.023,1463.195)(2762.042,1523.580)
	    (2707.385,1575.370)(2644.668,1617.037)(2575.742,1647.349)
	    (2502.643,1665.410)(2427.531,1670.688)(2352.625,1663.026)
	    (2280.138,1642.651)(2212.209,1610.164)(2150.847,1566.526)
	    (2097.864,1513.024)(2054.825,1451.241)(2023.000,1383.000)
	    \put(4302,915){\circle{16}}
	    \put(4307,911){\circle{1176}}
	    \drawline(2022,1387)(2022,443)
	    \drawline(2845,1387)(2845,443)
	    \drawline(1285,976)(1777,976)
	    \blacken\drawline(1711.340,955.485)(1777.000,976.000)(1711.340,996.515)(1711.340,955.485)
	    \drawline(3106,964)(3599,964)
	    \blacken\drawline(3533.340,943.485)(3599.000,964.000)(3533.340,984.515)(3533.340,943.485)
	    \drawline(12,1263)(12,1262)(13,1259)
	    (14,1255)(15,1248)(17,1239)
	    (20,1227)(24,1212)(29,1195)
	    (35,1176)(42,1155)(49,1132)
	    (58,1109)(69,1085)(81,1060)
	    (94,1035)(109,1009)(126,984)
	    (146,958)(168,932)(194,905)
	    (223,878)(257,850)(294,823)
	    (336,796)(382,770)(427,747)
	    (473,727)(519,709)(565,693)
	    (610,679)(655,667)(698,656)
	    (741,647)(782,640)(824,634)
	    (864,628)(904,624)(943,620)
	    (982,617)(1019,614)(1054,612)
	    (1088,611)(1119,610)(1147,609)
	    (1172,608)(1192,608)(1210,607)
	    (1223,607)(1232,607)(1238,607)
	    (1242,607)(1243,607)
	    \drawline(4088,703)(4087,702)(4084,699)
	    (4079,694)(4072,687)(4063,677)
	    (4052,664)(4039,649)(4024,631)
	    (4009,613)(3994,592)(3979,571)
	    (3965,549)(3953,526)(3942,502)
	    (3933,477)(3926,451)(3922,424)
	    (3922,395)(3925,365)(3933,334)
	    (3946,302)(3962,274)(3981,247)
	    (4000,223)(4020,201)(4038,181)
	    (4056,164)(4071,149)(4086,135)
	    (4099,124)(4111,113)(4123,103)
	    (4135,94)(4146,85)(4159,77)
	    (4173,68)(4189,60)(4208,51)
	    (4229,43)(4254,34)(4282,27)
	    (4314,20)(4350,15)(4388,12)
	    (4428,13)(4465,17)(4501,25)
	    (4535,34)(4566,44)(4594,55)
	    (4620,66)(4643,77)(4663,88)
	    (4681,98)(4698,108)(4713,117)
	    (4726,127)(4740,136)(4753,146)
	    (4765,156)(4778,167)(4792,180)
	    (4806,193)(4820,209)(4836,226)
	    (4851,245)(4867,267)(4882,291)
	    (4896,317)(4907,344)(4915,373)
	    (4918,404)(4917,435)(4911,464)
	    (4901,492)(4888,517)(4872,542)
	    (4854,565)(4833,587)(4812,608)
	    (4788,628)(4764,647)(4739,665)
	    (4715,683)(4690,699)(4667,714)
	    (4646,727)(4627,738)(4611,748)
	    (4598,755)(4578,767)
	    \blacken\drawline(4644.858,750.810)(4578.000,767.000)(4623.748,715.627)(4644.858,750.810)
	    \drawline(12,1263)(13,1263)(17,1264)
	    (23,1265)(32,1267)(44,1270)
	    (60,1273)(80,1276)(102,1279)
	    (127,1283)(155,1286)(184,1289)
	    (215,1291)(247,1291)(281,1291)
	    (316,1289)(352,1285)(391,1280)
	    (431,1271)(474,1260)(520,1246)
	    (568,1228)(618,1207)(669,1181)
	    (713,1156)(755,1130)(795,1103)
	    (832,1075)(868,1047)(901,1019)
	    (932,991)(961,963)(989,935)
	    (1015,907)(1040,879)(1063,851)
	    (1086,824)(1107,797)(1128,770)
	    (1147,745)(1164,720)(1181,698)
	    (1195,678)(1207,660)(1218,644)
	    (1227,632)(1233,622)(1238,615)
	    (1241,611)(1242,608)(1243,607)
	    \drawline(437,1014)(435,1014)(430,1014)
	    (421,1014)(408,1015)(391,1015)
	    (371,1015)(348,1014)(324,1013)
	    (299,1011)(273,1007)(248,1003)
	    (222,997)(197,989)(172,979)
	    (147,965)(124,949)(102,929)
	    (86,909)(72,889)(62,869)
	    (53,851)(47,835)(42,820)
	    (38,807)(35,795)(33,785)
	    (31,774)(30,764)(29,753)
	    (29,742)(30,730)(33,716)
	    (36,700)(42,684)(51,666)
	    (63,649)(78,632)(97,618)
	    (118,607)(139,598)(159,591)
	    (179,586)(196,583)(213,580)
	    (228,578)(242,577)(256,576)
	    (269,576)(284,576)(299,577)
	    (316,578)(335,580)(355,584)
	    (378,590)(402,597)(426,608)
	    (449,621)(472,639)(490,660)
	    (504,682)(514,704)(522,727)
	    (527,751)(529,775)(530,799)
	    (530,822)(529,845)(528,866)
	    (526,885)(524,901)(520,929)
	    \blacken\drawline(549.595,866.901)(520.000,929.000)(508.977,861.099)(549.595,866.901)
	    \drawline(2520,694)(2522,692)(2525,689)
	    (2531,682)(2541,673)(2554,661)
	    (2569,645)(2588,628)(2608,609)
	    (2630,590)(2653,570)(2678,550)
	    (2703,531)(2729,513)(2756,496)
	    (2785,480)(2815,466)(2847,454)
	    (2880,445)(2914,439)(2946,437)
	    (2976,438)(3002,442)(3025,446)
	    (3043,452)(3059,457)(3072,463)
	    (3083,469)(3092,474)(3100,481)
	    (3108,487)(3116,494)(3124,502)
	    (3133,512)(3143,523)(3152,536)
	    (3163,552)(3172,571)(3180,592)
	    (3185,615)(3186,639)(3183,663)
	    (3178,685)(3172,705)(3165,723)
	    (3158,738)(3152,751)(3145,763)
	    (3139,773)(3133,783)(3126,793)
	    (3118,803)(3109,814)(3098,825)
	    (3084,839)(3068,853)(3048,869)
	    (3025,886)(2999,903)(2969,918)
	    (2934,931)(2899,941)(2866,946)
	    (2833,949)(2802,949)(2772,946)
	    (2743,942)(2715,937)(2688,931)
	    (2663,924)(2639,917)(2619,910)
	    (2601,905)(2572,894)
	    \blacken\drawline(2626.116,936.468)(2572.000,894.000)(2640.668,898.105)(2626.116,936.468)
	    
\put(1283,621){\makebox(0,0)[lb]{\smash{{{\SetFigFont{10}{8.4}{\rmdefault}{\mddefault}{\updefault}$a$}}}}}
	    
\put(9,1311){\makebox(0,0)[lb]{\smash{{{\SetFigFont{10}{8.4}{\rmdefault}{\mddefault}{\updefault}$0$}}}}}
	    
\put(681,400){\makebox(0,0)[lb]{\smash{{{\SetFigFont{10}{8.4}{\rmdefault}{\mddefault}{\updefault}$f^n(l)$}}}}}
	    
\put(796,1160){\makebox(0,0)[lb]{\smash{{{\SetFigFont{10}{8.4}{\rmdefault}{\mddefault}{\updefault}$l$}}}}}
	    
\put(85,404){\makebox(0,0)[lb]{\smash{{{\SetFigFont{10}{8.4}{\rmdefault}{\mddefault}{\updefault}$R_f$}}}}}
	    
\put(4283,969){\makebox(0,0)[lb]{\smash{{{\SetFigFont{10}{8.4}{\rmdefault}{\mddefault}{\updefault}$0$}}}}}
	    
\put(3306,1063){\makebox(0,0)[lb]{\smash{{{\SetFigFont{10}{8.4}{\rmdefault}{\mddefault}{\updefault}$e$}}}}}
	    
\put(2298,1434){\makebox(0,0)[lb]{\smash{{{\SetFigFont{10}{8.4}{\rmdefault}{\mddefault}{\updefault}$\field{C} 
/ \field{Z}$ }}}}}
	    
\put(3000,650){\makebox(0,0)[lb]{\smash{{{\SetFigFont{10}{8.4}{\rmdefault}{\mddefault}{\updefault}$g$}}}}}   
\put(4350,151){\makebox(0,0)[lb]{\smash{{{\SetFigFont{10}{8.4}{\rmdefault}{\mddefault}{\updefault}$h$}}}}}
	  \end{picture}
	}
   \caption{Schematics of cylinder renormalization.}
  \end{center}
\end{figure}
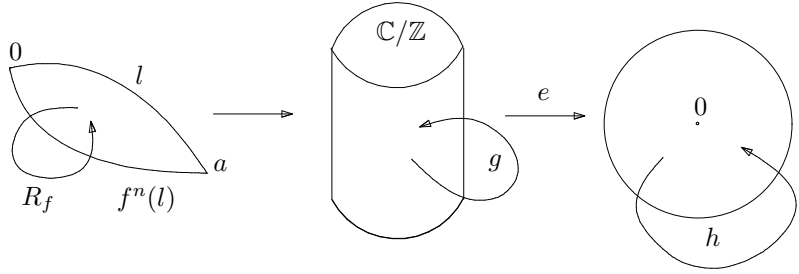

\noindent
{\bf Cylinder renormalization operator.}
The cylinder renormalization operator is defined as follows. 
Let $f\in{\aaa C}_W$.
 Suppose that for $n\in\NN$ there exists a simple arc $l$ which connects a fixed point 
$a$ of $f^n$ to $0$, and has the property that $f^n(l)$ is again a simple arc 
whose only intersection with $l$ is at the two endpoints. Let $C_f$ be the topological 
disk in $\CC\setminus\{0\}$ bounded by $l$ and $f^n(l)$.
We say that $C_f$ is a {\it fundamental crescent} if 
the iterate $f^{-n}|_{C_f}$ mapping $f^{n}(l)$ to $l$  is defined and univalent, and the quotient
of $\overline{C_f\cup f^{-n}(C_f)}\setminus \{0,a\}$ by the iterate $f^n$ is conformally
isomorphic to $\CC/\ZZ$. 
Let us denote $R_f$ the first return map of $C_f$, and let us denote  $z$ the 
critical point of this map (corresponding to the orbit of $0$). 
 Let $g$ be the map $R_f$ becomes under the above isomorphism, mapping $z$ to $0$,
and $h=e\circ g\circ e^{-1}$, where $e(z)=\exp{[-2 \pi i z]}$.
We say that $f$ is {\it cylinder renormalizable with period $n$}, if $h\in {\aaa C}_V$ for 
some $V$, and call $h$ a {\it cylinder renormalization} of $f$ (see Figure 1).

We summarize below the basic properties of cylinder renormalization proven in \cite{Ya1}:

\begin{prop}
\label{properties cyl ren}
Suppose $f\in{\aaa C}_W$ is cylinder renormalizable, and its renormalization $h_f$ is contained
in ${\aaa C}_V$. Denote $C_f$ the fundamental crescent corresponding to the
renormalization. Then the following holds.
\begin{itemize}
\item Every other fundamental crescent $C'_f$ with the same endpoints as $C_f$, and such that
$C'_f\cup C_f$ is a topological disk, produces the same renormalized map $h_f$. 
\item There exists an open neighborhood $U(f)\subset{\aaa C}_W$ such that
every map $g\in U(f)$ is cylinder renormalizable, with a fundamental crescent $C_g$ 
which can be chosen to move continuously with $g$.
\item Moreover,  the dependence $g\mapsto h_g$ of the
cylinder renormalization on the map $g$ is an analytic mapping ${\aaa C}_W\to{\aaa C}_V$.
\end{itemize}
\end{prop}

We now want to discuss the dynamical properties of the cylinder renormalization 
of maps with Siegel disks
derived in \cite{Ya3}. To simplify the exposition let us specialize to the case
when the rotation number of the Siegel disk is the golden mean $\theta=(\sqrt{5}+1)/2$.
The golden mean  is represented by an infinite continued fraction 
$$\theta=1+\cfrac{1}{1+\cfrac{1}{1+\cfrac{1}{\cdots}}}\equiv 1+[1,1,1,\ldots].$$
As is customary, we will denote $p_n/q_n$ the $n$-th convergent
$$p_n/q_n=\underbrace{[1,1,1,\ldots,1]}_n.$$

\begin{thm}[\cite{Ya3}]
\label{thm-renorm}
There exists a space ${\bf C}_U$ and an analytic mapping $\hat f\in {\bf C}_U$ 
which has a Siegel disk $\Delta_{\theta}$
with rotation number $\theta$ whose boundary is a quasicircle passing through the
critical point $1$, such that the following holds:
\begin{itemize}
\item[(I)] There exists a branch of cylinder renormalization with period $q_k$, $k\in\NN$,
which we denote $\cren$ such that
$$\cren\hat f=\hat f;$$
\item[(II)] the quadratic polynomial  $P_\theta(z)=e^{2\pi i\theta}z+z^2$
is infinitely cylinder renormalizable, and 
$$\cren^k P_\theta\to \hat f,$$
at a uniform geometric rate;
\item[(III)] 
the cylinder renormalization $\cren$ is an analytic and compact operator mapping a neighborhood
of the fixed point $\hat f$ in ${\aaa C}_{U}$ to ${\aaa C}_{U}$. Its linearization $\cL$ at 
$\hat f$ is a compact operator, with at least one 
 eigenvalue with the absolute value greater than one.
\end{itemize}
\end{thm}

A central open questions in the study of $\cren$ is the following:

\begin{conj}
Except for the one unstable eigenvalue, the rest of the spectrum of $\cL$ is compactly contained
in the unit disk.
\end{conj}

Our numerical study of $\cren$ will begin with 
empirically establishing the convergence to $\hat f$. We will then 
make  explicit the choice of the neighborhood $U$
in the above Theorem. Experimental evidence suggests that it can be taken as a round disk
$\DD_\rho$ for some particular value of $\rho$. Having numerically established this, we will
then proceed to experimentally verify the Conjecture.

\section{Construction of the conformal isomorphism to the cylinder} \label{Cylinder}

The principal difficulty in numerical, as well as analytical, study of cylinder renormalization is the
non-explicit nature of the conformal isomorphism 
$$\Phi:\overline{C_f\cup f^{-n}(C_f)}\setminus\{0,a\}\underset{\approx}{\longrightarrow}\CC^*$$
of a fundamental crescent, which is a part of the definition of $\cren$.
An analytic approach to this construction based on the Measurable Riemann Mapping Theorem was 
presented by the first author in \cite{Gay}. It has its roots in the complex-dynamical folklore;
similar arguments are found, for instance, in the work of Lyubich \cite{Lyu} and Shishikura \cite{Shish}. 

In \cite{Gay}, the first author demonstrates how this approach can be implemented constructively,
with rigorous error bounds. We will give a brief outline here.

\medskip
\noindent
{\bf Uniformization of the cylinder using the Measurable Riemann Mapping Theorem}.
We will start with a description of our choice of a fundamental crescent $C_f^n$ with period $q_n$ 
for a map  $f\in{\bf C}_U$ sufficiently close to $\hat f$.

To construct the boundary curve $l_n$ of $C_f^n$ consider first the union $\tl l_n$ of 
two  parabolas $x +i (Ax^2+B x)$ and $(C y^2+D y +E)+i y$: the first passing through points $0$ and $f^{q_{n+2}+q_n}(1)$, the second --- through   $f^{q_{n+2}+q_n}(1)$ and a repelling fixed point $a_{q_n}$. All parameters in these two parabolas are defined uniquely after one specifies their common tangent line at   $f^{q_{n+2}+q_n}(1)$ (see equations $(\ref{l_1})$ below). While somewhat arbitrary, this choice has the virtue of possessing 
a simple analytic form. It can be shown rigorously (see \cite{Gay}) that by modifying $\tl l_n$ in sufficiently small neighborhoods
of the endpoints (small enough not to influence our numerical experiments) we obtain a curve $l_n$ which 
together with $f^{-q_n}(l_n)$ bounds a  fundamental crescent $C_f^n$ for $\cren$.

Now consider the following conformal change of coordinates for $z\in C_f^n$:
\begin{equation}\label{tau}
  z=\tau(\xi)={a_{q_n} \over 1-e^{i \alpha  \xi+\beta }}, \quad \tau^{-1}(z)={1 \over i a}\left[\ln{\left(1-{a_{q_n} \over z } \right) }-\beta \right].
\end{equation}

The normalizing constant $\beta$ will be chosen so that 
 $$\tau^{-1}(f^{q_{n+2}}(1))=0,$$ 
while a real positive $\alpha$ will be specified by the condition
 $$|\tau^{-1}(f^{q_{n+2}+q_n}(1))|=1\text.$$ 
The choice of of this coordinate is motivated by the fact that $\tau^{-1}$ maps the interior of the 
fundamental crescent $C_f^n$ conformally onto the interior of an infinite vertical closed strip
 $\mathcal{S}$, whose width 
is comparable to one (cf. Figure 2). 
Next, similarly to \cite{Shish},  define a function 
$$\tilde{g}_n:\mathcal{U}\equiv\left\{u+i v \in \field{C}\right.: 0 \le \Re {w} \le 1\left. \!\!\!\!\!\!\phantom{C}\right\}\longrightarrow\mathcal{S}$$
 by setting
\begin{align}
  \nonumber   \tilde{g}_n(u+i v)&= (1-u) \tau^{-1} (f^{-q_n}(\gamma_n(v)))+u \tau^{-1} (\gamma_n(v)),
\end{align}
where $\gamma_n$ is a parametrization 
$$\gamma_n : \field{R}  \to  l_n$$
which we will specify below. 

Let  $\sigma_0$  be the standard conformal structure on $\field{C}$, and let $\sigma=\tilde{g}_n^* \sigma_0$ be its pull-back on $\mathcal{U}$. Extend this conformal structure to $\field{C}$ through 
$$\sigma\equiv (T^k)^* \sigma\text{ on  }T^{-k}(\mathcal{U}),\text{ where  }T(w)=w+1\text{, for all }k\in\NN.$$
Assuming the mapping $\tl g_n$ is quasiconformal, the dilatation of $\sigma$ is bounded in the plane.
By Measurable Riemann Mapping Theorem (see e.g. \cite{AB}) there exists a unique quasiconformal 
mapping $\tilde{g}: \field{C} \mapsto  \field{C}$ such that $\tilde{g}^* \sigma_0 = \sigma$, normalized so that $\tilde{g}(0)=0$ and $\tilde{g}(1)=1$. Notice that $\tilde{g} \circ T \circ \tilde{g}^{-1}$  preserves the standard conformal structure:
\begin{equation}
\nonumber \left(\tilde{g} \circ T \circ \tilde{g}^{-1}\right)^* \sigma_0 = (\tilde{g}^{-1})^* 
\circ T^* \circ \tilde{g}^* \sigma_0 =(\tilde{g}^{-1})^* \circ T^* \sigma=(\tilde{g}^{-1})^* \sigma =(\tilde{g}^*)^{-1} \sigma=\sigma_0,
\end{equation}
and therefore it is a conformal automorphism of $\field{C}$. Liouville's Theorem implies that this mapping is affine.
By construction, it does not have any fixed points in $\CC$, and hence, is a translation.
Finally, $\tilde{g} \circ T \circ \tilde{g}^{-1}(0)=1$, and thus  $$\tilde{g} \circ T \circ \tilde{g}^{-1}\equiv T.$$

By the definition of $\tilde{g}_n$, 
$$\tilde{g}_n^{-1} \circ \tau^{-1} \circ f^{q_n} \circ \tau=T \circ \tilde{g_n}^{-1}$$
 on the image of $f^{-q_n}(l_n)$ by $\tau^{-1}$. 
Set $\phi=\tilde{g} \circ \tilde{g}_n^{-1}$, and $\tilde\Phi\equiv\phi \circ \tau^{-1}$.
Clearly, $\tilde\Phi \mod \ZZ$ is a
 desired conformal isomorphism  
$$\overline{C_{f}^n\cup f^{-q_n}(C_f^n)}\setminus\{0,a_{q_n}\}\underset{\approx}{\longrightarrow}\field{C} / \field{Z}.$$

Again, set $e(z)=e^{-2\pi iz}$, and  $g=e \circ \tilde{g} \circ e^{-1}$. Since 
\begin{equation}
\nonumber {g_{\bar{z}}(e(w)) \over g_z(e(w))}={e(w) \over \overline{e(w)}} {\tilde{g}_{\bar{w}}(w) \over \tilde{g}_w(w)},  
\end{equation}
the $1$-periodic function $\tilde{g}$ is a solution of the Beltrami equation 
\begin{equation}
\nonumber \tilde{g}_{\bar{w}}=\tilde{\mu}\tilde{g}_{w}, \quad \tilde{\mu}=(\tilde{g}_n)_{\bar{w}}/(\tilde{g}_n)_{w}
\end{equation}
whenever $g$ is a solution of  
\begin{equation}
g_{\bar{z}}=\mu g_{z}, \quad \mu(z)=(z/\bar{z}) \tilde{\mu}(e^{-1}(z)).
\end{equation}

\noindent
Thus, we have reduced the problem of finding 
$$\Phi\equiv e \circ \tilde\Phi=g \circ e \circ \tilde{g}^{-1}_n \circ \tau^{-1}$$
 to that of finding the properly normalized solution of the Beltrami equation
\begin{equation}\label{beltt}
g_{\bar{z}}=\mu g_z, \quad \mu(z)={z \over \bar{z}} {(\tilde{g}_n)_{\bar{w}} (e^{-1}(z)) \over (\tilde{g}_n)_{w} (e^{-1}(z)) } 
\end{equation}
on the punctured plane $\field{C}^*$. 

It remains to describe the choice of the parametrization of $l_n$ in the definition of $\tilde g_n$.
It is convenient for us to parametrize $l_n$ using the radial coordinate in $\CC$.
For $n=1$ and $f \in {\bf C}_U$, sufficiently close the empirical fixed point of the cylinder renormalization with period $1$
we use the following parametrization:
\begin{equation}\label{l_1}
\lambda_1(r) = \left\{ (x(r),A x(r)^2+B x(r)), \quad r \le \tilde{r},  ~ \atop  (C y(r)^2+D y(r)+E,y(r)) , \quad r>\tilde{r},  \right. 
\end{equation}
where 
\begin{eqnarray}
 \nonumber x(r)&=&{\Re{f^4(1)} \over |f^4(1)  |} T(r), \\  
  \nonumber y(r)&=&\Im{f^4(1)} {|a_1-f^4(1)|+|f^4(1)|-T(r) \over |a_1-f^4(1)|}+\Im{a_1} {T(r)-|f^4(1)|) \over |a_1-f^4(1)|},\\
  \nonumber T(r)&=& {|a_1-f^4(1)|+|f^4(1)| \over \sqrt{r} +1} \sqrt{r},
 \end{eqnarray}
and $\tilde{r}$ is defined through $T(\tilde{r})=|f^4(1)|$. Constants $A$, $B$, $C$, $D$ and $E$ are fixed by the conditions $0, f^4(1), a_1 \in  l_1$, together with the requirement that the slope of the  common tangent line to both parabolas at the point  $f^4(1)$ is equal to $1.1$.

 This particular choice of the parametrization is  motivated  by the speed of 
convergence of the iterative scheme in the Measurable Riemann Mapping Theorem.

We define the following function on $\field{C}^*$:
\begin{equation} \label{g_1}
 g_1(r,\phi)= \left(\eta(-\phi)+{\phi \over 2 \pi}\right) \tau^{-1}(f^{-1}(\lambda_1 (r)))+\left( 1-\eta(-\phi)- {\phi \over 2 \pi} \right)\tau^{-1}(\lambda_1 (r)),
\end{equation}
where $-\pi < \phi \le \pi$, and $\eta$ is the Heaviside step function (we have adopted the convention $\eta(0)=1$). Then, according to $(\ref{beltt})$, the Beltrami differential $\mu$ is given by the following expression:
\begin{equation}\label{belt2}
\mu(r e^{i\phi})=e^{2 i \phi} {r \partial_r g_1(r,\phi)+i \partial_{\phi} g_1(r,\phi) \over r \partial_r g_1(r,\phi)-i \partial_{\phi} g_1(r,\phi) } 
\end{equation}
on $\field{C} \setminus (-\infty,0]$. This is the expression that we have used to compute the Beltrami differential in our numerical studies.

\medskip
\noindent
{\bf A constructive Measurable Riemann Mapping Theorem (MRMT).}
To solve the Beltrami equation numerically, we use the constructive MRMT proved by the first author in \cite{Gay}.
Before formulating it, we need to recall two integral operators used in the 
 classical approach to the proof of MRMT (see \cite{AB}).

The first of them is  Hilbert Transform:
\begin{equation}\label{Hilbert_tr}
 T [h](z)={i \over 2 \pi} \lim_{\epsilon \to 0} \int \int_{\field{C} \setminus B(z,\epsilon)} {h(\xi) \over (\xi -z)^2} \ d \bar\xi \wedge d \xi, 
\end{equation}
the second is Cauchy Transform
\begin{equation}\label{Cauchy_tr}
 P [h] (z)={i \over 2 \pi} \int \int_{\field{C}} {h(\xi) \over (\xi -z)} \ d \bar\xi \wedge d 
\xi.
\end{equation}

Hilbert Transform is a well-defined  bounded operator on 
$L_p(\field{C})$ for all $2<p<\infty$. For every such $p$ there exists a constant 
$c_p$ such that the following holds (cf. \cite{CalZyg}):
$$\left\| T [h] \right\|_p \le c_p \left\| h \right\|_p\text{ for any }h \in L_p(\field{C})\text{, and }
c_p \rightarrow 1\text{ as }p \rightarrow 2.$$ 

We are now ready to state the constructive Measurable Riemann Mapping Theorem of \cite{Gay}:
\begin{thm}\label{constructive MRMT}
Let $\mu \in L_\infty(\bar{\field{C}})$ and an integer $p>2$  be such that $\left\| \mu \right\|_\infty\le K < 1$ and $K c_p<1$, where 
\begin{equation}
\nonumber c_p= \cot^2(\pi/2p).
\end{equation}
 Assume that $\mu=\nu+\eta+\gamma$, where $\nu$ and $\eta$ are compactly supported in $\field{D}_R$, and  $\gamma(z)$ is supported in 
$\bar{\field{C}} \setminus \field{D}_R$.  Furthermore, let $\eta$ be in $L_p(\field{D}_R)$ and $\left\| \eta \right\|_p < \delta$ for 
some 
sufficiently small $\delta$. Also, let $h^* \in L_p(\field{C})$ and $\epsilon$   be such that $B_p(h^*,\epsilon)$, the ball of radius 
$\epsilon$ around $h^*$ in  $L_p(\field{C})$, contains  $B_p(T[\nu (h^*+1) ], c_p \epsilon')$, with 
\begin{equation}
\nonumber \epsilon'= \delta \  {\rm essup}_{\field{D}_R}{ |h^*  +1|}+K  \epsilon.
\end{equation}
Then the solution $g^\mu$ of the Beltrami equation $g^\mu_{\bar{z}}=\mu g^\mu_z$ admits the following bound:
\begin{equation}\label{mainbound}
\left|g^\mu(z)-g^\nu_*(z)\right|\le F(\epsilon',R; z,g^\nu_*(z),p,K,c_p)
\end{equation}
where $g_*^\nu(z)=P[\nu(h^*+1)](z)+z$, and $F(\epsilon',R)=O(\epsilon',R^{-4/p})$ is an explicit function of its arguments.
\end{thm}

Given the theorem, the algorithm for producing an approximate solution of the Beltrami equation is as follows.
Given  a $\mu$ as in the condition of the theorem, we first iterate 
\begin{equation}\label{Hilbert_iter}
h \rightarrow T [ \nu (h+1)]
\end{equation}
 to find a numerical approximation $h^*_{a}$ to the solution of the equation $T [ \mu 
(h^*+1)]=h^*$. After that, we compute an approximate solution as  
\begin{equation}\label{solution}
g^\nu_{a}(z)=P [\nu (h^*_{a}+1)](z)+z.
\end{equation}

One can obtain rigorous computer-assisted bounds on such solution using Theorem $\ref{constructive MRMT}$. Such bounds have been indeed implemented  in  \cite{Gay} for a particular case of the golden mean quadratic polynomial. However, in the present numerical work we will not require such estimates.

In the Appendix, we will discuss several numerical algorithms for 
the two integral transforms appearing in this scheme.

\medskip\section{Empirical convergence to a fixed point}\label{Results}

\begin{figure}

\centerline{\resizebox{100mm}{!}{\includegraphics{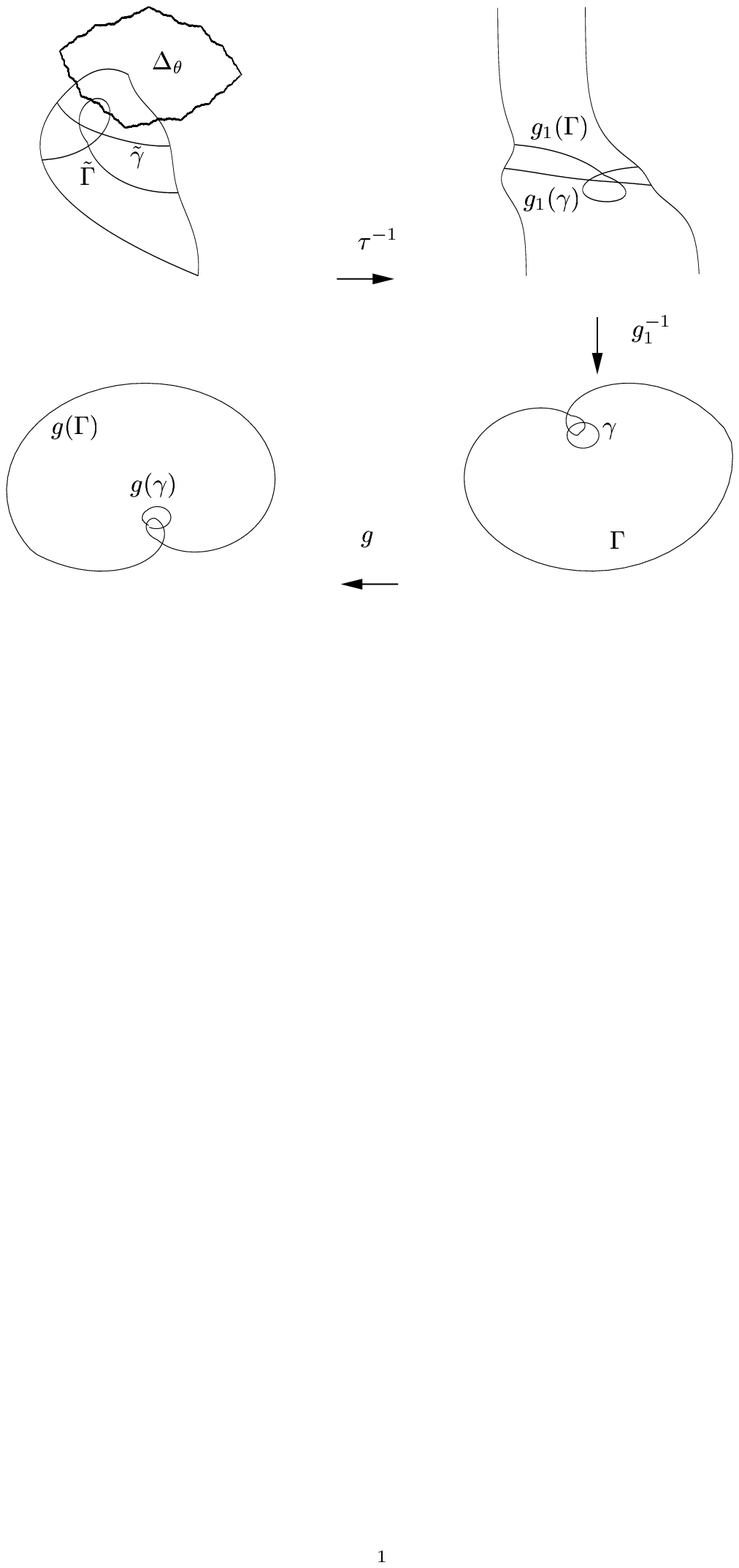}} }
\caption{Schematics of renormalization.
The contours $g(\gamma)\mapsto g(\Gamma)$ are used to find a polynomial approximation of  $\cren f$ 
through the Cauchy integral formula.}     
\label{Contours} 
\end{figure}

An appropriate choice of the domains of analyticity for the renormalized functions is central to a successful numerical
 implementation of
 cylinder renormalization. Our numerical approximation to the renormalization fixed point is a finite-degree truncation 
of a function analytic 
in $\DD_3$ (see Section \ref{Domain} for a detailed explanation of this choice of the domain). 
However, for the purposes of obtaining bounds 
on higher-order terms, we will consider a smaller analyticity domain, a disk of radius $\rho=2.266$. Thus the 
cylinder renormalization will be {\it a priori} an analytic operator in a neighborhood of the fixed point in the 
Banach space ${\bf L}^1_\rho$
with $\rho=2.266$.
\footnote{We have implemented the procedure for the cylinder renormalization described in Section \ref{Cylinder} and a
 particular method of solving the Beltrami equation (see Appendix) as a set of routines in  in the programming language 
Ada 95 (cf \cite{ADA} for the language standard). We have parallelized our programs and compiled them with the public version 
3.15p of the GNAT compiler \cite{GNAT}. The programs (\cite{programs})  have been run on the computational cluster of 92 2.2 
GHz AMD Opteron processors located at the University of Texas at Austin.}


Given an $f \in {\bf L}^1_\rho$, a numerical approximation to its cylinder renormalization of order $n$ is built as follows.  
As the first step, we construct a
fundamental crescent $C_{f}^n$ as described above,  and find the normalized solution $g$ of the Beltrami equation 
$$g_{\bar{z}}(z)=\mu(z) g_z(z)\text{ with }\mu\text{ as in }(\ref{belt2})$$ 
as described in Section \ref{Cylinder} and the Appendix. 
Next, we choose a contour $\gamma$ in the domain of $g$, 
and map this contour into the fundamental crescent by $\tau \circ g_1$:  
$$\tilde{\gamma}\equiv\tau \circ g_1(\gamma).$$
 Applying the first
return map to the points of this contour, we 
obtain $\tilde{\Gamma}=R_{f}(\tilde{\gamma})$, and find the images,  $g(\gamma)$ and $g(\Gamma)$ (where $\Gamma=g^{-1}_1 
(\tau^{-1}(\tilde{\Gamma}))$). 
The coefficients in a finite order polynomial approximation to 
$$\cren f= g \circ g^{-1}_1\circ \tau^{-1}  \circ R_{f} \circ \tau \circ g_1 \circ g^{-1} $$
are then found via the Cauchy integral formula using these two contours  $g(\gamma)$ and  $g(\Gamma)$ (see Figure \ref{Contours}).

As seen in \thmref{thm-renorm}, the sequence of the cylinder renormalizations of the quadratic polynomial  $P_{\theta}$ converges to a
 fixed point $\hat f=\cren \hat f$. We have used this fact to compute an approximate renormalization fixed point 
$\hat f_a$ as the  cylinder renormalization $\cren^k P_{\theta}$ of order  $k=11$.

Further, we improved this approximation by iterating 
\begin{equation}
\label{iter}
\hat f_a \mapsto \field{P}_s \circ  \cren \hat f_a 
\end{equation}
 where $\field{P}_s$ is the projection on the candidate stable manifold of $\hat f$
$${\mathrm W}^s=\{f \in {\bf L}^1_\rho: f'(0)=e^{2 \pi \theta i} \},$$
defined by setting
$$\field{P}_s[f](x) \equiv f(x)+(e^{2 \pi \theta i}-f'(0))x.$$
  
  In this way we have obtained a polynomial $\hat f_a$ of degree $17$, and have estimated that, not taking into account the errors in the
solution of the Beltrami equation and due to round-off,
\begin{equation}\label{error}
\|\cren \hat f_a-\hat f_a\|_\rho \le 1.88 \times 10^{-3} \approx 0.89\times 10^{-4}||\hat f_a||_\rho.
\end{equation}
Moreover,
the iteration (\ref{iter}) does not lead to a significant variation in the computed 
values for the coefficients of $\hat f_a$, which indicates that the original approximation
is indeed quite accurate. The largest change is in the highest coefficient, which 
differs by $0.4\%$ for $\hat f_a$ and its renormalization. Of course, this represents a
negligible correction to the  absolute
value of the coefficient itself.

The approximate expression for $\hat f_a$ is as follows (all numbers truncated to show six significant digits):
\begin{eqnarray*}
\hat f_a(x)=&x&e^{2\pi i\theta}+\\	
 &x^2&(\phantom{-}8.00882 \times 10^{-1}+ i4.07682 \times 10^{-1})+\\
 &x^3&(-4.12708 \times 10^{-1} +i2.97670 \times 10^{-2})+\\
&x^4& (\phantom{-} 1.02033 \times 10^{-1} -i9.83702 \times 10^{-2})+\\
 &x^5&(\phantom{-} 2.61573 \times 10^{-5} +i4.13871 \times 10^{-2})+\\
&x^6&(-8.42868 \times 10^{-3} -i6.96474 \times 10^{-3})+\\
 &x^7&(\phantom{-} 2.60095 \times 10^{-3} -i6.58544 \times 10^{-4})+\\
&x^8&(-2.01382 \times 10^{-4} +i5.95113 \times 10^{-4})+\\
&x^9&(-9.40057 \times 10^{-5} -i1.11237 \times 10^{-4})+\\
&x^{10}&( \phantom{-}3.21762 \times 10^{-5} -i4.40144 \times 10^{-6})+\cdots
\end{eqnarray*}

\section{Domain of Analyticity of the Renormalization Fixed Point} \label{Domain}

\noindent
\begin{figure}[h]
  \begin{center}
 \resizebox{80mm}{!}{\includegraphics{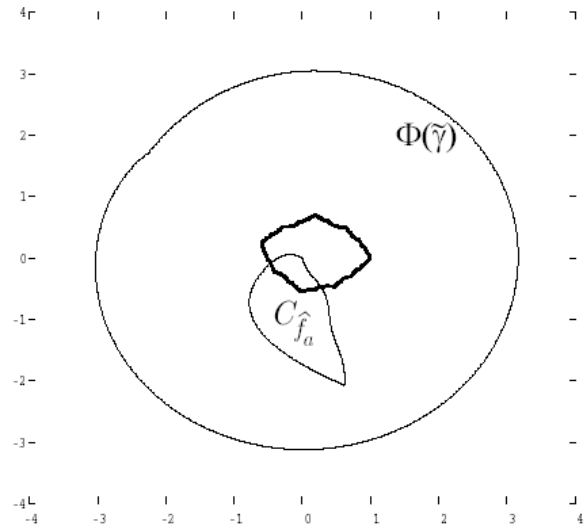}} 
  \caption{ The fundamental domain $C_{\hat{f}_a}$ together with $\Phi(\tilde{\gamma})$.}   \label{Domain-fig}
  \end{center}
\end{figure}

\begin{figure}
  \begin{center}
 \resizebox{105mm}{!}{\includegraphics{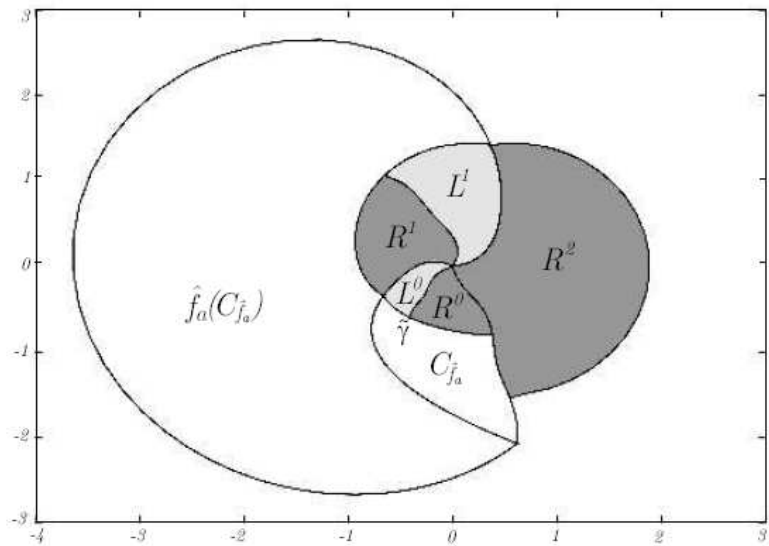}} 

\setlength{\unitlength}{0.00083333in}
%
\caption{The orbit of $C^0_{\hat{f}_a}$: the orbit of $L^0_{\hat{f}_a}$, 
$L^k\!\!\equiv{\hat{f}_a}^k(L^0_{\hat{f}_a})$, is rendered in light gray, that of 
$R^0_{\hat{f}_a}$, 
$R^k\!\!\equiv{\hat{f}_a}^k(R^0_{\hat{f}_a})$, --- in dark gray. 
}   \label{Orbit}
  \end{center}
\end{figure}

\medskip
\noindent
{\bf Compactness of $\cren$.} 
We have verified experimentally the compactness property of $\cren$ 
stated in \thmref{thm-renorm}. More precisely, we observe the following empirical fact:

\medskip
\noindent
{\it Set $\rho =2.266$ and $\rho'=3$. Then we can take $U\equiv \DD_\rho$ in \thmref{thm-renorm}. 
More specifically, the fixed point $\hat f$ is a well-defined analytic mapping in ${\bf C}_\rho$,
and moreover, if we denote 
$$g\equiv \hat f|_{\DD_\rho},\text{ then }\cren g\in {\bf C}_{\rho'}.$$}


To verify the claim numerically, we have used the approximation $\hat f_a$ obtained in
the previous section. To estimate $\rho'$, we have chosen a curve
$\tilde{\gamma}$  in 
the fundamental crescent, such that $\Phi_{\hat{f}_a}(\tilde{\gamma})$ is a simple closed loop 
that encircles $\field{D}_{3}$ (Figure 3). We then verify that the orbit under the return map of the 
component $C^0_{\hat{f}_a}$ of the set  $C _{\hat{f}_a} \setminus \tilde{\gamma}$, such that $0 
\in 
\partial C^0_{\hat{f}_a}$, lies within  $\field{D}_{2.266}$ (see Figure \ref{Orbit}). For 
our 
choice of the curve $\tilde{\gamma}$, the return map of the set  
$$C^0_{\hat{f}_a}=L^0_{\hat{f}_a} 
\cup 
R^0_{\hat{f}_a}$$ is given by the $2$-nd and the $3$-d iterates of $\hat{f}_a$ on 
$L^0_{\hat{f}_a}$ 
and $R^0_{\hat{f}_a}$, respectively.

\section{Hyperbolic properties of cylinder renormalization}\label{Hyperbolicity}

\noindent
{\bf The expanding direction of} $\cren$.
It is not difficult to see that the operator $\cren$ possesses an expanding direction at $\hat f$ (cf. \cite{Ya3}):

\begin{proof}[Proof of \thmref{thm-renorm}, Part (III)]
Let $v(z)$ be a vector field in ${\aaa C}_U$,
$$v(z)=v'(0)z+o(z).$$
Denote $\gamma_v$ the quantity 
$$\gamma_v=\frac{v'(0)}{\hat f'(0)}={e^{-2\pi i\theta}}v'(0).$$
For a smooth family
$$\hat f_t(z)=\hat f(z)+tv(z)+o(t),$$
we have 
$$\hat f_t(z)=\alpha^v_t(z)(\hat f'(0)z+o(z)),\text{ where }\alpha_t(0)=1+t\gamma_v+o(t).$$
The $q_{m+1}$-st iterate
$$\hat f_t^{q_{m+1}}(z)=(\alpha^v_t(z))^{q_{m+1}}((\hat f'(0))^{q_{m+1}}z+o(z)).$$
In the neighborhood of $0$ the renormalized vector field $\cL v$ is obtained
by applying a uniformizing coordinate 
$$\Phi(z)=(z+o(z))^\beta,\text{ where }\beta=\frac{1}{\theta q_m\mod 1}.$$
Hence,
$$\alpha_t^{\cL v}(0)=[(\alpha_t^v(0))^{q_{m+1}}]^\beta,$$
so
$$\gamma_{\cL v}=\Lambda \gamma_v,\text{ where }\Lambda=\beta q_{m+1}>1.$$
Hence the spectral radius
$$R_{\text{Sp}}(\cL v)>1,$$
and since every 
 non-zero element of the spectrum of a compact operator is
an eigenvalue, the claim follows.
\end{proof}

\medskip
\noindent
{\bf Numerical verification of hyperbolicity of $\cren$.}
It is natural to conjecture:

\begin{conj}\label{Conj}
There exists an open neighborhood $\cU\subset {\bf C}_U$ containing $\hat f$ such that
$\cren$ is a strong contraction in $$W=\{f\in \cU\;|\;f'(0)=e^{2\pi i\theta}\}.$$
Thus, $W=W^s_{\text{loc}}(\hat f)$.
\end{conj}

To verify this conjecture numerically, we have to justify using a finite-dimensional
approximation to $\cL$ to test for contraction. For this we rely on a numerical
observation discussed in the previous section:
$$\cL:{\bf L}^1_\rho\to{\bf L}^1_{\rho'},\text{ with }\rho=2.266,\text{ and }\rho'=3.$$
This implies, that the finite-dimensional approximations of $\cL$ obtained by truncating all
the powers higher that $z^N$ will converge geometrically fast in $N$.

Set $h_j$ to be the coordinate vectors $h_j(z)=z^j / \rho^j$, so that
$$||\cL ||_\rho=\sup||\cL h_j||_\rho.$$
Since a perturbation $\hat f+\eps h_j$ does not lie in ${\bf L}^1_\rho$, we perturb along a different set of vectors:
\begin{equation}
e_j={g_j \over \| g_j \|_\rho}, \quad  g_j(z)=z^j-{j \over j+1} z^{j+1}, \quad  j \ge 1,
\end{equation}
which form a basis in ${\bf L}^1_\rho$.

Numerically, to estimate the spectral radius $$R_{\text{Sp}}\left(\cL\arrowvert_{T_{\hat f}W}\right)$$ we can fix a large enough $N$
and a small $\eps$ (we have used the value $\eps=0.01$), compute for each $e_j$, $2\leq j \leq  N$, the finite difference
$$\frac{1}{\eps}(\cren (\hat f_a+\eps e_j)-\cren\hat f_a),$$
truncate past the $N$-th power -- and expand over the basis vectors $e_j$ to obtain an $(N-1)\times (N-1)$
matrix $A_N$. Below we present the approximate expression for $A_6$ (the numbers have been truncated to
the fifth decimal).
\medskip
\begin{center}
\begin{tabular}{|l|l|l|l|l|}
\hline
\phantom{$-$}$0.45879-$   & \phantom{$-$}$0.68789-$   &\phantom{$-$}$0.11338-$  &\phantom{$-$}$0.13041+$ &$0.15824+$\\
\phantom{$-$}$0.97624i$   & \phantom{$-$}$0.46254i$   &\phantom{$-$}$0.09738i$  &\phantom{$-$}$0.07490i$ &$0.11616i$\\
\hline
            $-0.13666+$   &             $-0.64474+$   &\phantom{$-$}$0.33937-$  &\phantom{$-$}$0.14710-$ &$-0.21006$+\\
\phantom{$-$}$1.72834i$   &\phantom{$-$}$0.54306i $   &\phantom{$-$}$0.50837i $ &\phantom{$-$}$0.13849i $&\phantom{$-$}$0.02552i$ \\
\hline
$-0.90634-$&\phantom{$-$}$0.27155-$&$-0.05765+$&$-0.22948+$&\phantom{$-$}$0.18338-$\\
\phantom{$-$}$1.37322i $&\phantom{$-$}$0.38270i $&\phantom{$-$}$1.09081i $&\phantom{$-$}$0.14700i $&\phantom{$-$}$0.39078i$ \\
\hline
\phantom{$-$}$1.23970+$&$-0.08549+$&$-0.68227-$&\phantom{$-$}$0.12817-$&\phantom{$-$}$0.10981+$\\
\phantom{$-$}$0.20634i $&\phantom{$-$}$0.23219i $&\phantom{$-$}$0.81153i $&\phantom{$-$}$0.14685i $&\phantom{$-$}$0.47861i $\\
 \hline
 $-0.63443+$&$-0.02489-$&\phantom{$-$}$0.80205+$&$-0.04014+$&$-0.34685-$ \\
\phantom{$-$}$0.58168i $&\phantom{$-$}$0.18893i $&\phantom{$-$}$0.00357i $&\phantom{$-$}$0.12864i  $&\phantom{$-$}$0.19936i $       \\
\hline
\end{tabular}
\end{center}

\medskip

\begin{center}
Table 1. Matrix $A_6$.
\end{center}

\medskip

This matrix has the spectral radius
$$R_{\text{Sp}}(A_6)\approx 0.53.$$

\noindent
{\bf Estimating the spectral radius.} We now proceed to produce a justification for  the above numerical experiment.
We will equip ${\bf L}^1_\rho$, viewed as a vector space,  with a new $l_1$-norm
\begin{equation}
| f |_\rho=\sum^{\infty}_{k=1} |f_k|,
\end{equation}
where $f_k$ are the coefficients in the expansion of $f$ in the basis $\{ e_j \}$: $f=\sum^{\infty}_{k=1} f_k e_k$;
and denote the new Banach space by ${\bf \tl{L}}^1_\rho$. The projection $\field{P}_{\le N}$ on ${\rm span}_{1 \le  j \le N} 
\{e_j\}$ will be defined by setting 
\begin{equation} \label{P_N}
\field{P}_{\le N} f= \sum^{N}_{j=1}f_j e_j.
\end{equation}
We will also abbreviate $\field{I}-\field{P}_{\le N}$ as $\field{P}_{>N}$.

We would like to emphasize that $\field{P}_{\le N} f \in {\bf {\tilde L}}^1_{\rho}$ whenever $f \in {\bf {\tilde L}}^1_{\rho}$, 
and therefore the operator
\begin{equation}\label{A}
\cA=\field{P}_{\le N} \cL \field{P}_{\le N}
\end{equation}
serves as a finite-dimensional approximation to the action of $\cL$ on ${\bf \tilde L}^1_\rho$. We will now make the latter 
statement more precise.

To this end, notice, that 
\begin{equation}\label{A+H}
\cL=\cA +\cL \field{P}_{> N}  +\field{P}_{> N} \cL \field{P}_{\le N} =\cA+\cH.
\end{equation}

The following Lemma demonstrates how one can obtain an upper bound on the spectral radius of the differential $\cL$ at
the fixed point in terms of the norm of a power of the finite-rank  operator $\cA$ and the magnitude of the norm of 
$\cH$. 

\begin{lemma}\label{Spectral_radius}
Let $\cL=\cA+\cH$ be a bounded operator on some Banach space, such that $\| \cA^k \| < \gamma <1$ for some $k \ge 
1$ and, $\| \cH \| <\delta <1$.   Then, the spectral radius $R_{\text{Sp}}(\cL)$ satisfies
\begin{equation} \label{R_sp}
R_{\text{Sp}}(\cL) \le \gamma^{1/k} (1+C \delta / \gamma)^{1/k}
\end{equation}
for some (explicit) constant $C$.
\begin{proof}
The claim follows from the spectral radius formula. First,
\begin{equation}
\nonumber R_{\text{Sp}}(\cL)=\overline{\lim}_{n \rightarrow \infty} \| \cL^n \|^{1 \over n} =\overline{\lim}_{n \rightarrow \infty} \left\| \cL^{k \left[{ n \over k} \right]+k \left\{ {n \over k} \right\} } \right\|^{1 \over n}\! \le\overline{\lim}_{n \rightarrow \infty} \left\| \cL^{k \left[{ n \over k} \right] }\right\|^{1 \over n}\overline{\lim}_{n \rightarrow \infty} \| \cL^k \|^{1-{ k \left[{n \over k} \right] \over n}  }.
\end{equation}
The norm $\| \cL^k \|$ is finite, and therefore
\begin{equation}
\nonumber \overline{\lim}_{n \rightarrow \infty} \| \cL^k \|^{n- k \left[{n \over k} \right] \over n}=1.
\end{equation}
Then
\begin{equation}
\nonumber R_{\text{Sp}}(\cL) \le \overline{\lim}_{n \rightarrow \infty} \left\| \cL^{k \left[{ n \over k} \right] }\right\|^{1 \over n} \!\!  \le \overline{\lim}_{n \rightarrow \infty} \left\| \cL^{k \left[{ n \over k} \right] }\right\|^{1 \over \left[  {n \over k} \right] k }  \overline{\lim}_{n \rightarrow \infty}  \| \cL^k \|^{{ \left[  {n \over k} \right] k -n } \over n k } \le  \overline{\lim}_{m \rightarrow \infty} \left\| \cL^{k m }\right\|^{1 \over m k }.
\end{equation}
 Let $\Prefix_{n}{C}_k$ denotes the binomial coefficients, and let $C=\sum^{k}_{i=1} \Prefix_{k} C_i \| \cA ^{k-i} \| \|\cH\|^{i-1} $ then
\begin{eqnarray}
\nonumber R_{\text{Sp}}(\cL) &\le&    \overline{\lim}_{n \rightarrow \infty} \left[ \|\cA^k \|^n +\sum^{n}_{i=1} \Prefix_{n} C_i \|\cA^k\|^{n-i} (C 
\delta)^i \right]^{1 \over k n}\\
\nonumber  & \le & \overline{\lim}_{n \rightarrow \infty} \left[ \gamma^n (1 + ((C \delta/\gamma+1)^n-1) \right]^{1 \over k n} \\
\nonumber &= & \gamma^{1/k} (1+C \delta / \gamma)^{1/k}.
\end{eqnarray} 
\end{proof} 
\end{lemma}

It is left now to bound the difference of $\cL$ from $\cA$. First, we state the  following Cauchy-type estimate, whose straightforward proof will be left to the reader:

\begin{prop}\label{Derivative}
Assume that an operator $\cren$ is analytic in an open ball $B_r(\hat f)$ $\subset {\bf \tilde L}^1_\rho$.
 Let $\epsilon<1$, and $h \in{\bf \tilde L}^1_\rho $ be such that $| h |_\rho <r$. Then
\begin{equation} \label{approx_der}
| \cren (\hat f+\epsilon h)- \cren \hat f-\epsilon \cL  h  |_\rho \le 
{\epsilon^2 \over 1-\epsilon}  \sup_{|s|\le 1}  | \cren ( \hat f+s h)-\cren \hat f |_\rho.
\end{equation}
\end{prop}

Note that
\begin{equation} \label{der_norm}
\nonumber |\cL \arrowvert_{T_{\hat f} W} |_\rho \le \sup_{j \ge 2} | \cL  e_j   |_\rho.
\end{equation}
This, together with the  preceding Proposition and the compactness property of renormalization, immediately implies that  $| \cL \arrowvert_{T_{\hat f} W} |_\rho$  can be bound by a finite 
difference. Specifically, for all $j > N$:
\begin{eqnarray} \label{C}
\nonumber | \cL e_j |_{\rho}  &  \le & \epsilon^{-1} \sup_{h \in \field{P}_{> N} {\bf \tilde{L}}^1_\rho ~\atop |h|_\rho \le 1}   | \cren(\hat f +\epsilon h)- 
\hat f |_{\rho}+  {\epsilon \over  1-\epsilon} \sup_{h \in \field{P}_{> N} {\bf \tilde{L}}^1_\rho ~\atop |h|_\rho \le 1}   |\cren (\hat f+h)- \hat f |_{\rho}   \\
 \nonumber  &  \le & \epsilon^{-1}  \sup_{h \in \field{P}_{> N} {\bf \tilde{L}}^1_\rho ~\atop |h|_\rho \le 1}  
 \!\!|\field{P}_{\le N} [ \cren(\hat f +\epsilon h)- \hat f] |_{\rho} + {\epsilon \over 1-\epsilon} \sup_{h \in \field{P}_{> N} {\bf \tilde{L}}^1_\rho ~\atop |h|_\rho \le 1} \!\! |\field{P}_{\le N} [\cren (\hat f+h)- \hat f] |_{\rho} \\
\nonumber  &  + & \!\left( {\rho \over \rho'} \right)^{N+1} \!\!\! \!\!\!\!\! \epsilon^{-1}\!\! \sup_{h \in \field{P}_{> N} {\bf \tilde{L}}^1_\rho ~\atop |h|_\rho \le 1}  \!\!\! |\field{P}_{> N} [ \cren(\hat f \!+\!\epsilon h)\! -\! \hat f] |_{\rho'}\!+\! {\epsilon \over  1-\epsilon} \sup_{h \in \field{P}_{> N} {\bf \tilde{L}}^1_\rho ~\atop |h|_\rho \le 1} \!\!\! |\field{P}_{> N} [\cren (\hat f\!+\!h)\!-\! \hat f] |_{\rho'}  \\
&\equiv& C_1,
\end{eqnarray}
Similarly, for all $2 \le j\le N$:
\begin{eqnarray} \label{Cp}
\nonumber | \field{P}_{> N} \cL e_j |_{\rho} \!   & \le & \! \epsilon^{-1} \sup_{h \in {T_{\hat f} W} ~\atop |h|_\rho \le 1} \!  
 |\field{P}_{> N}  [ \cren(\hat f +\epsilon h)- \hat f] |_{\rho} + {\epsilon \over 1-\epsilon }  \sup_{h \in {T_{\hat f} W} ~\atop |h|_\rho \le 1} | \field{P}_{> N} 
[\cren (\hat f+h)- \hat f] |_{\rho}  \\
\nonumber    & \le & \!\! \left( {\rho \over \rho'} \right)^{N+1} \!\!\! \!\!\!\!\! \epsilon^{-1}\!\!\!\sup_{h \in {T_{\hat f} W} ~\atop |h|_\rho \le 1} \!\!  |\field{P}_{> N}[\cren(\hat f \!+\! \epsilon h)\! -\! \hat f] |_{\rho'} \! + \! {\epsilon \over 1-\epsilon }  \sup_{h \in {T_{\hat f} W} ~\atop |h|_\rho \le 1}\!\! | \field{P}_{> N} [\cren (\hat f\!+\!h)\!-\! \hat f] |_{\rho'}\\
& \equiv & C_2.
 \end{eqnarray}

We would like to emphasize that these bounds use the fact that $\cL$ is a compact operator in an essential way.  Bounds $(\ref{C})$ and $(\ref{Cp})$  can be used to estimate $| (\cL-\cA) \arrowvert_{T_{\hat f} W}|_\rho$. In particular, according to equation $(\ref{A+H})$,
\begin{equation}
\nonumber | (\cL-\cA) \arrowvert_{T_{\hat f} W } |_\rho  \le  \max \left\{ | \cL \field{P}_{> N} |_\rho , |\field{P}_{> N} \cL  
\arrowvert_{T_{\hat f} W}  \field{P}_{\le N} |_\rho \right\} \leq \max \left\{ C_1,C_2 \right\}.
\end{equation}
This expression provides a bound on $\delta$ in $(\ref{R_sp})$.

We have chosen $N=14$, and experimentally bounded $ C_1$ and $C_2$  by testing on  vectors $h=e_{15}$, and $h=e_{2}$ which empirically maximize the respective suprema. Choosing $k=80$ in $(\ref{R_sp})$, we have the following values for the constants that enter 
 estimate $(\ref{R_sp})$:
\begin{equation}
\nonumber  \gamma < 2.07 \times 10^{-18}, \quad \delta < 0.24, \quad C < 8.4 \times 10^{-6},
\end{equation}
therefore, according to $(\ref{R_sp})$, 
\begin{equation}
\nonumber R_{\text{Sp}}(\cL \arrowvert_ {T_{\hat f} W}) < 0.85.
\end{equation}

A better bound can be obtained if one computes all relevant constants for a larger value of $N$, which requires more computer time. 
It is plausible that the spectral radius is close to $0.58$, since we have observed that as $N$ increases the largest eigenvalue of the operator 
$$\field{P}_{\le N} \cL \arrowvert_{T_{\hat f} W}  \field{P}_{\le N}$$
 converges to  
\begin{equation}
\nonumber \lambda=0.15+i 0.56.
\end{equation}
This eigenvalue has been truncated to $2$ decimal places.

As a final comment, note that the following simple observation implies that perturbations in the directions of the
vectors $h_j$ can also be used for estimating the spectral radius of $\cL\arrowvert_{T_{\hat f}W}$. 

\begin{prop}
We have
$$\text{Spec}(\cL\arrowvert_{{\bf B}^1_\rho}) = \text{Spec}(\cL\arrowvert_{{\bf 
L}^1_\rho}).$$
\end{prop}

\noindent
To see this, note that the only difference between the spectra is that $0$ (contained in 
both spectra) is an eigenvalue  of the operator $\cL\arrowvert_{{\bf B}^1_\rho}$ corresponding to linear rescalings.  
We leave the straightforward details to the reader.

\medskip\section{Appendix}\label{Appendix}

The objective of this Appendix will be to describe how the Cauchy (\ref{Cauchy_tr}) and Hilbert (\ref{Hilbert_tr}) transforms can be computed numerically. 

The Constructive Measurable Riemann Theorem \ref{constructive MRMT} deals with $L_p$ functions which generally do not need to be differentiable. Therefore, one has to chose an appropriate representation of the $L_p$ functions that enter the Theorem \ref{constructive MRMT}; possibly, as a collection of values on a grid, or as a Fourier series with the radially dependent coefficients. The latter choice has been made, for instance, in \cite{Daripa3}, \cite{Gay}, \cite{GayKhmel}, and will be also adopted in the present paper.     

 Represent $h$ and $P [h]$ in (\ref{Cauchy_tr}) as:
\begin{eqnarray}
\label{h_coeff} h(r e^{i \theta})&=&\sum^{\infty}_{k=-\infty} h_k(r) e^{i k \theta}, \\
\label{Pt_coeff} P [ h ](r e^{i \theta})&=&\sum^{\infty}_{k=-\infty} p_k(r) e^{i k \theta},
\end{eqnarray}
where the coefficients of the $P$-transform are given by 
\begin{equation} \label{pk}
p_k(r)= {1 \over 2 \pi } \int^{2 \pi}_0 e^{-i k \theta} P [h](r e^{i \theta}) \ d \theta.
\end{equation}

A classical theorem of analysis (cf \cite{Ahlfors}, \cite{Markovic}) states that Cauchy 
transform of an $L_p$-function, $p > 2$, is well-defined and is H\"older continuous with 
exponent $1-2/p$. In \cite{Daripa3} and \cite{GayKhmel} this fact has been used to show that 
the Fourier coefficients of Cauchy transform are given by the following equations 
\begin{equation}\label{Cauchy_alg}
 p_k(r)=\left\{
\begin{array}{cc}
  \displaystyle    2\int_0^r\left(\frac{r}{\rho}\right)^k h_{k+1}(\rho)d\rho,&k < 0,\\
  \displaystyle   -2\int_r^\infty\left(\frac{r}{\rho}\right)^k h_{k+1}(\rho)d\rho,&k\ge0.
\end{array}
 \right.
\end{equation}

To obtain similar formulae for Hilbert transform, assume that $h$ is a H\"older continuous 
function compactly supported in an open  disk around zero of radius $R$, $B(0,R) \subset \field{C}$.  The Hilbert transform of such function is known to exist as a Cauchy principal value (cf \cite{Ahlfors}, \cite{Carleson}). As with Cauchy transform, represent this transform as a Fourier series:
\begin{eqnarray}\label{T_coeff} 
T [ h ](r e^{i \theta})&=&\sum^{\infty}_{k=-\infty} c_k(r) e^{i k \theta}, \quad c_k(r)= {1 \over 2 \pi } \int^{2 \pi}_0 e^{-i k \theta} T[h](r e^{i \theta}) \ d \theta. 
\end{eqnarray}

In \cite{Daripa1}, \cite{Daripa2} and \cite{GayKhmel} the authors arrive at the following expressions for these coefficients:
\begin{eqnarray}
\label{Hilbert_alg_1}c_0(0)&=& -2 \lim_{\epsilon \rightarrow 0} \int^R_\epsilon {h_2(\rho) \over \rho}  d \rho, \ {\rm and} \ c_k(0)=0, \ {\rm whenever} \ k \ne 0, \\
\label{Hilbert_alg_2}c_k(r)&=&A_k \int^r_0 { r^k  \over \rho^{k+1} }  h_{k+2}(\rho)   d \rho+B_k\int^R_r  { r^k  \over \rho^{k+1} }  h_{k+2}(\rho)   d \rho +h_{k+2}(r),
\end{eqnarray}
where 
\begin{eqnarray} \label{A_k_B_k} 
A_k&=&\left\{ 0 \ , \ k \ge 0,  \atop  2 (k+1) \ , k<0, \right.   \ {\rm and} \    B_k=\left\{-2  (k+1),  \ k \ge 0, \atop 0, \ k<0. \right.
\end{eqnarray}

We would like to mention that the fact that Hilbert transform is a singular integral operator makes a rigorous justification of the formulas (\ref{Hilbert_alg_1})--(\ref{Hilbert_alg_2}) significantly more involved than that of (\ref{Cauchy_alg}).

Formulas (\ref{Cauchy_alg}) and (\ref{Hilbert_alg_1})--(\ref{Hilbert_alg_2}) can be used to construct an efficient algorithm for solving a Beltrami equation. Given  values of $h$, for instance, on a circular $N \times M$ grid that contains the compact support of $h$, one can use a fast Fourier transform (FFT, cf \cite{NR}) to find the values of the coefficients $h_k$ at the radii $r_i$, $ 1 \le i\le M$. Next, one can use these values to construct a piecewise constant, a piecewise linear or a spline approximation of the functions $h_k$ (the choice of approximation, of course, depends on the known or expected smoothness of $h_k$). This allows one to compute integrals in (\ref{Cauchy_alg}) and in (\ref{Hilbert_alg_1})--(\ref{Hilbert_alg_2}). Armed with these implementations of Hilbert and Cauchy transforms, one can try to solve the Beltrami equation (\ref{beltt}), first by running iterations (\ref{Hilbert_iter}) for some time, and, finally, applying (\ref{solution}). It is convenient to use the point-wise multiplication of grid values of $h$ and $\mu+1$ inside Hilbert transform in (\ref{Hilbert_iter}), rather than the multiplication of their Fourier series: the order of the computational complexity of the point-wise multiplication is $O(NM)$, as opposed to $O(NM^2)$ for the series. The transition from the representation of $h$ as a Fourier series to point values at each iteration step can be performed with the help of the  FFT. This way, the computational complexity of one iteration step becomes $O(N M \log_2 M)$.


\begin{thebibliography}{****} 
\bibitem[Ahl]{Ahlfors}  L. Ahlfors, {\it Lectures on quasiconformal mappings}, Van Nostrand-Reinhold, Princeton, New Jersey, 1966.

\bibitem[AB]{AB} L. Ahlfors, L. Bers, Riemann's mapping theorem for variable metrics, {\it  Ann. of Math. } (2) {\bf 72}(1960), 385--404.


\bibitem[Ber]{Bers}  L. Bers, Quasiconformal mappings, with applications to differential equations, function theory and topology, {\it  American Mathematical Society}  {\bf 83}(1977), 1083--1100.

\bibitem[BI]{Boyar_Iwan} B. V. Boyarskii and T. Iwanier, Quasiconformal mappings and nonlinear elliptic equations in two variables, I and II {\it  Bull. Acad. Polon. Sci., S\'er. Math. Astronom. Phys.} {\bf 12}(1974),  473--478 and 479--484.

\bibitem[Bo]{Boyarskii} B. V. Boyarskii, Generalized solutions of systems of differential equations of first order and elliptic type with discontinuous coefficients, {\it  Mat. Sb. N. S.}  {\bf 43(85)}(1957),  451--503.

\bibitem[CZ]{CalZyg} A. P. Calderon, A. Zygmund,  On singular integrals, {\it  Amer. J. Math.}  {\bf 78}(1956), 289--309.

\bibitem[Car]{Carleson}  L. Carleson, Th. W. Gamelin, {\it Complex Dynamics}, Springer (1991).

\bibitem[Da1]{Daripa3}   P. Daripa, A fast algorithm to solve non-homogeneous Cauchy-Riemann equations in the complex plane, {\it   SIAM J. Sci. Statist. Comput.} {\bf 13}(1992) 1418--1432. 

\bibitem[Da2]{Daripa1}   P. Daripa, A fast algorithm to solve the Beltrami equation with applications to quasiconformal mappings,  {\it  J. Comput. Phys.}  {\bf 106}(1993) 355--365. 

\bibitem[Da3]{Daripa2}   P. Daripa and D. Mashat, Singular Integral Transforms and Fast Numerical Algorithms, {\it   Numer. Algor.}  {\bf 18}(1998) 133--157. 

\bibitem[dF1]{dF1} E. de Faria, {\it Proof of universality for critical circle mappings}, Thesis, CUNY, 1992.

\bibitem[dF2]{dF2} E. de Faria, Asymptotic rigidity of scaling ratios for critical circle mappings, {\it Ergodic Theory Dynam. Systems} {\bf 19}(1999), no. 4, 995--1035.

\bibitem[Gai]{Gay}  D. Gaydashev,  Computer-Assisted Bounds on the Solution of a Beltrami Equation and Applications to Renormalization, e-print math.DS/0510472 at Arxiv.org.

\bibitem[GK]{GayKhmel}D. Gaydashev, D. Khmelev, On Numerical Algorithms for the Solution of a Beltrami Equation, e-print math.DS/0510516 at Arxiv.org.

\bibitem[Lan1]{La1} O.E. Lanford, Renormalization group methods for critical circle mappings with general rotation number, {\it VIIIth International Congress on Mathematical Physics (Marseille,1986)}, World Sci. Publishing, Singapore, 532--536, (1987).

\bibitem[Lan2]{La2} O.E. Lanford, Renormalization group methods for critical circle mappings. Nonlinear evolution and chaotic phenomena, {\it NATO adv. Sci. Inst. Ser. B:Phys.} {\bf 176}(1988), Plenum, New York, 25--36.

\bibitem[Lyu]{Lyu} M. Lyubich, Dynamics of rational transformations: topological picture, {\it  Uspekhi Mat. Nauk}  {\bf 41}(1986),  no. 4(250), 35--95

\bibitem[Mar]{Markovic}  V. Markovic, {\it Quasiconformal maps}, Lecture notes taken by A. Fletcher. 
\begin{verbatim}http://www.maths.warwick.ac.uk/~fletcher/qcmaps.pdf \end{verbatim}

\bibitem[MN]{MN} N.S. Manton, M. Nauenberg, Universal scaling behaviour for iterated maps in
the complex plane, {\it Commun. Math. Phys.} {\bf 89}(1983), 555--570.

\bibitem[MP]{MP} R.S. MacKay, I.C. Persival, Universal small-scale structure near the boundary of Siegel disks of arbitrary rotation numer, {\it Physica} {\bf 26D}(1987), 193--202.

\bibitem[McM]{McM} C. McMullen, Self-similarity of Siegel disks and Hausdorff dimension of  Julia sets, {\it Acta Math.} {\bf 180}(1998), 247-292.


\bibitem[NR]{NR}   W. H. Press, B. P. Flannery, S. A. Teukolsky  and   W. T. Vetterling, {\it Numerical Recipes in Fortran. The Art of Scientific Computing}, Cambridge: Cambridge University Press 1992.

\bibitem[Shi]{Shish} M. Shishikura, The Hausdorff dimension of the boundary of the Mandelrot set and Julia sets, {\it   Ann. of Math.} {\bf 43}(1942), 607--612.

\bibitem[Sieg]{Sieg} C.L. Siegel, Iteration of analytic functions, {\it Ann. Math.} {\bf 43}(1942), 607-612.

\bibitem[Stir]{Stir} A. Stirnemann,  Existence of the Siegel disc renormalization fixed point,  {\it  Nonlinearity } {\bf 7}(1994),  no. 3, 959--974.

\bibitem[Wi]{Wi} M. Widom, Renormalisation group analysis of quasi-periodicity in analytic maps,
{\it Commun. Math. Phys.} {\bf 92}(1983), 121-136.

\bibitem[Ya1]{Ya1} M. Yampolsky, Hyperbolicity of renormalization of critical circle maps,
{\it  Publ. Math. Inst. Hautes Etudes Sci.} {\bf 96}(2002), 1--41.

\bibitem[Ya2]{Ya2} M. Yampolsky, Renormalization horseshoe for critical circle maps,
{\it  Commun. Math. Physics} {\bf 240}(2003), 75--96.

\bibitem[Ya3]{Ya3} M. Yampolsky,  Siegel disks and renormalization fixed points, e-print math.DS/0602678 at Arxiv.org

\bibitem[ADA1]{ADA}   S. T. J. Taft and R. A. Duff (eds), {\it Ada 95 Reference Manual: Language and Standard Libraries, International Standard ISO/IEC 8652:1995(E)},  Lec. Notes in Comp. Science {\bf 1246}. 

\bibitem[ADA2]{GNAT}   Ada Core Technologies, 73 Fifth Ave, New York, NY 10003, USA.\\ See also {\tt ftp://cs.nyu.edu/pub/gnat}.

\bibitem[Prog]{programs} {\tt http://www.math.toronto.edu/gaidash/Programs/siegel-numerics.tar.bz2}.
\end{thebibliography}
\end{document}